\documentclass[12pt]{amsart}
\usepackage{latexsym,amssymb,amscd}
\def\frk{\frak}               % font for "Fraktur"

\def\mm{{\frk m}}

\def\B'c{{\mathcal{B'}}}
\def\U'c{{\mathcal{U'}}}

\def\xb{{\bold x}}
\def\yb{{\bold y}}
\def\zb{{\bold z}}
\def\opn#1#2{\def#1{\operatorname{#2}}} % to make operators
\def\bigin{{\rm bigin}}
\def\char{{\rm char}}
\opn\chara{char}
\opn\length{\ell}
\opn\projdim{proj\,dim}
\opn\injdim{inj\,dim}
\opn\ini{in}
\opn\rank{rank}
\opn\Tiefe{Tiefe}
\opn\grade{grade}
\opn\depth{depth}
\opn\height{height}
\opn\embdim{emb\,dim}
\opn\codim{codim}

\opn\Tr{Tr}
\opn\bigrank{big\,rank}
\opn\superheight{superheight}\opn\lcm{lcm}
\opn\trdeg{tr\,deg}%
\opn\reg{reg}
\opn\lreg{lreg}
\opn\hdeg{hdeg}
\opn\ideg{ideg}

\opn\div{div}
\opn\Div{Div}
\opn\cl{cl}
\opn\Cl{Cl}
\opn\Spec{Spec}
\opn\Supp{Supp}
\opn\supp{supp}
\opn\Sing{Sing}
\opn\Ass{Ass}
\opn\Ann{Ann}
\opn\Rad{Rad}
\opn\Soc{Soc}
\opn\Ker{Ker}
\opn\Coker{Coker}
\opn\Im{Im}
\opn\Hom{Hom}
\opn\Tor{Tor}
\opn\Ext{Ext}
\opn\End{End}
\opn\Aut{Aut}
\opn\id{id}

\opn\nat{nat}
\opn\GL{GL}
\opn\SL{SL}
\opn\mod{mod}
\opn\ord{ord}
\opn\Gin{Gin}
\opn\aff{aff}
\opn\con{conv}
\opn\relint{relint}
\opn\st{st}
\opn\lk{lk}
\opn\cn{cn}
\opn\core{core}
\opn\vol{vol}
\opn\gr{gr}

\let\union=\cup

\let\dirsum=\oplus
\let\tensor=\otimes
\let\iso=\cong

\let\Dirsum=\bigoplus

\def\pnt{{\raise 0.5mm\hbox{\large\bf.}}}
\def\lpnt{{\hbox{\large\bf.}}}
\let\to=\rightarrow

\def\Implies{\ifmmode\Longrightarrow \else
     \unskip${}\Longrightarrow{}$\ignorespaces\fi}
\def\implies{\ifmmode\Rightarrow \else
     \unskip${}\Rightarrow{}$\ignorespaces\fi}
\def\iff{\ifmmode\Longleftrightarrow \else
     \unskip${}\Longleftrightarrow{}$\ignorespaces\fi}

\let\:=\colon
\newtheorem{Theorem}{Theorem}[section]
\newtheorem{Lemma}[Theorem]{Lemma}
\newtheorem{Corollary}[Theorem]{Corollary}
\newtheorem{Proposition}[Theorem]{Proposition}
\newtheorem{Remark}[Theorem]{Remark}

\let\epsilon=\varepsilon
\let\phi=\varphi
\let\kappa=\varkappa
\title{Homological properties of bigraded algebras}
\author{Tim R\"omer}

\begin{document}

\maketitle

\begin{abstract}
We investigate the $x$- and $y$-regularity of a bigraded
$K$-algebra $R$ as introduced in \cite{ARCRNE}.
These notions are used to study asymptotic properties 
of certain finitely generated bigraded modules. As
an application we get for any equigenerated graded ideal $I$
upper bounds for the number $j_0$
for which $\reg(I^j)$ is a linear function for $j \geq j_0$.  
Finally, 
we give upper bounds for the $x$- and 
$y$-regularity of generalized Veronese algebras.
\end{abstract}

\section*{Introduction}
Let $S=K[x_1,\ldots,x_n,y_1,\ldots,y_m]$ be 
a standard bigraded 
polynomial ring with 
$\deg(x_i)=(1,0)$ and $\deg(y_j)=(0,1)$, and let 
$J \subset S$ be a bigraded ideal. 
In this paper we study homological properties of the
bigraded algebra $R=S/J$. 

First we consider the $x$- and the $y$-regularity of $R$.
According to \cite{ARCRNE} they are defined as follows:
$$
\reg_x^S(R)=\max\{a \in \mathbb{Z} \colon \beta_{i,(a+i,b)}^S(R)\neq 0 \text{ for some } i,b \in \mathbb{Z}  \},  
$$
$$
\reg_y^S(R)=\max\{b \in \mathbb{Z}  \colon \beta_{i,(a,b+i)}^S(R)\neq 0 \text{ for some } i, a \in \mathbb{Z}  \}
$$
where $\beta_{i,(a,b)}^S(R)=\dim_K \Tor_i^S(K,R)_{(a,b)}$ is
the $i^{\rm th}$ bigraded Betti number of $R$ in bidegree $(a,b)$.
We give a homological characterization of these regularities 
similarly as in the 
graded case (see \cite{ARHE2}). As an application 
we generalize a result of Trung \cite{TR} concerning $d$-sequences.
Furthermore we prove that 
$$
\reg_x^S(S/J)=\reg_x^S(S/\bigin(J))
$$ 
where $\bigin(J)$ is the bigeneric initial ideal of $J$ with respect to
the bigraded reverse lexicographic order induced 
by $y_1>\ldots>y_m>x_1>\ldots>x_n$.

It was shown in \cite{CUHETR} (or \cite{KO}) that for $j \gg 0$,
$\reg(I^j)$ is a linear function $cj+d$ in $j$ for a graded ideal
$I$ in the polynomial ring.
In \cite{KO} the constant $c$ is described in terms of invariants of $I$.
In this paper we give, in case $I$ is equigenerated, 
bounds $j_0$ such that for $j \geq j_0$
the function is linear and give also
a bound for $d$. Our methods 
can also be applied to $\reg(S^j(I))$, where
$S^j(I)$ is the $j^{\rm th}$ symmetric power of $I$.
 
In the last section we introduce a generalized Veronese algebra in the bigraded
setting. For a bigraded algebra $R$ and $\Tilde{\Delta}=(s,t) \in \mathbb{N}^2$ with 
$(s,t) \neq (0,0)$ we set
$$
R_{\Tilde{\Delta}}=\Dirsum_{(a,b)\in\mathbb{N}^2} R_{(as,bt)}.
$$
In the same manner as it is done for diagonal subalgebras in
\cite{COHETRVA}, we prove
that for these algebras  
$$
\reg_x^{S_{\Tilde{\Delta}}}(R_{\Tilde{\Delta}})=0 \text{ and }
\reg_y^{S_{\Tilde{\Delta}}}(R_{\Tilde{\Delta}})=0,
\text{ if }  s \gg 0 \text{ and } t \gg 0.
$$

\section{Preliminaries}
Throughout this paper, let $K$ be an infinite field and 
$S=K[x_1,\ldots,x_n,y_1,\ldots,y_m]$ be 
a standard bigraded 
polynomial ring with 
$\deg(x_i)=(1,0)$ and $\deg(y_j)=(0,1)$.
Let $M$ be a finitely generated bigraded $S$-module. 
For some bihomogeneous $w \in M$ with $\deg(w)=(a,b)$
we set $\deg_x(w)=a$ and $\deg_y(w)=b$.
Sometimes we will consider the $\mathbb{Z}$-graded modules 
$M_{(a,\ast)}=\Dirsum_{b \in \mathbb{Z} } M_{(a,b)}$ 
or $M_{(\ast,b)}=\Dirsum_{a \in \mathbb{Z} }M_{(a,b)}$.
If in addition $M$ is $\mathbb{N}^n \times\mathbb{N}^m$-graded, we write $M_{(u,v)}$
for the homogeneous component in bidegree $(u,v)$ where $u \in \mathbb{N}^n$ and $v \in \mathbb{N}^m$.
For $u \in \mathbb{N}^n$ we set $\supp(u)=\{i \colon u_i>0 \}$.

Define $\mm_x=(x_1,\ldots,x_n)=(\xb)$, $\mm_y=(y_1,\ldots,y_m)=(\yb)$ and $\mm=\mm_x+\mm_y$. 
Let $S_x=K[x_1,\ldots,x_n]$ and $S_y=K[y_1,\ldots,y_m]$ be the polynomial rings
with respect to the $x$-variables and the $y$-variables.

For some $u=(u_1,\ldots,u_n) \in \mathbb{N}^n$ and $v=(v_1,\ldots,v_m) \in \mathbb{N}^m$
we write $x^u y^v$ for the monomial 
$x_1^{u_1}\ldots x_n^{u_n}y_1^{v_1}\ldots y_m^{v_m}$. 
For $u, u' \in \mathbb{N}^n$ let $u \preceq u'$, if $u_i \leq u_i'$ for all $i$.
Furthermore we set $|u|=u_1+\ldots+u_n$. 
Let $\epsilon_i=(0,\ldots,0,1,0,\ldots,0)\in \mathbb{N}^n$ where the
entry $1$ is at the $i^{\rm th}$ position. For $t \in \mathbb{N}$ define $[t]=\{1,\ldots,t\}$.

We consider bigraded algebras $R=S/J$, which are quotients of $S$ by some 
bigraded ideal $J$. For a finitely generated bigraded $R$-module $M$ and $a,b \in \mathbb{N}$
let $\beta_{i,(a,b)}^R(M)=\dim_K \Tor^R_i(M,K)_{(a,b)}$ 
be the $i^{\rm th}$ bigraded Betti number in bidegree $(a,b)$.
We recall from \cite{ARCRNE} that
$$
\reg_x^R(M)=\sup\{a \in \mathbb{Z} \colon \beta_{i,(a+i,b)}^R(M)\neq 0 \text{ for some } i,b \in \mathbb{Z}  \},  
$$
$$
\reg_y^R(M)=\sup\{b \in \mathbb{Z} \colon \beta_{i,(a,b+i)}^R(M)\neq 0 \text{ for some } i,a \in \mathbb{Z} \}
$$
is the $x$- and $y$-regularity of $M$. 
For $R=S$ we set $\reg_x(M)=\reg_x^S(M)$ and $\reg_y(M)=\reg_y^S(M)$. 

Let 
$K_{\lpnt}(k,l;M)$ denote the Koszul complex of $M$
and
$H_{\lpnt}(k,l;M)$ the Koszul homology of $M$
with respect to $x_1,\ldots,x_k$ and $y_1,\ldots,y_l$ (see \cite{BRHE} for details).
If it is clear from the context, we write $K_{\lpnt}(k,l)$ and $H_{\lpnt}(k,l)$ instead of
$K_{\lpnt}(k,l;M)$ and
$H_{\lpnt}(k,l;M)$.
Note that $K_{\lpnt}(k,l;M)=K_{\lpnt}(k,l;S)\tensor_S M$
where $K_{\lpnt}(k,l;S)$ is the exterior algebra on
$e_1,\ldots,e_k$ and $f_1,\ldots,f_l$ with $\deg(e_i)=(1,0)$ and $\deg(f_j)=(0,1)$
together with a differential $\partial$
induced by $\partial(e_i)=x_i$ and $\partial(f_j)=y_j$.
For a cycle $z \in K_{\lpnt}(k,l;M)$ we denote with $[z]\in H_{\lpnt}(k,l;M)$
the corresponding homology class.
There are two long exact sequences relating the homology groups:
$$
\ldots \to H_i(k,l;M)(-1,0)\overset{x_{k+1}}{\to} H_i(k,l;M)\to H_i(k+1,l;M) \to H_{i-1}(k,l;M)(-1,0) 
$$
$$
\overset{x_{k+1}}{\to} \ldots \to
H_0(k,l;M)(-1,0)\overset{x_{k+1}}{\to} H_0(k,l;M)\to H_0(k+1,l;M)\to  0
$$
and
$$
\ldots \to H_i(k,l;M)(0,-1)\overset{y_{l+1}}{\to} H_i(k,l;M)\to H_i(k,l+1;M) \to H_{i-1}(k,l;M)(0,-1) 
$$
$$
\overset{y_{l+1}}{\to} \ldots \to
H_0(k,l;M)(0,-1)\overset{y_{l+1}}{\to} H_0(k,l;M)\to H_0(k,l+1;M)\to 0.
$$
The map $H_i(k,l;M)\to H_i(k+1,l;M)$ is induced by the inclusion of
the corresponding Koszul complexes. Every homogeneous element $z \in K_{\lpnt}(k+1,l;M)$
can be uniquely written as $e_{k+1}\wedge z'+z''$ with $z',z''\in K_{\lpnt}(k,l;M)$. Then
$H_i(k+1,l;M) \to H_{i-1}(k,l;M)(-1,0)$ is given 
by sending $[z]$ to $[z']$.
Furthermore
$H_i(k,l;M)(-1,0)\overset{x_{k+1}}{\to} H_i(k,l;M)$ is just the multiplication
with $x_{k+1}$. The maps in the other exact sequence are analogue.
  
\section{Regularity}
Let $R$ be a bigraded algebra.
To simplify the notation we do not distinguish between the 
polynomial ring variables $x_i$ or $y_j$ and the
corresponding residue classes in $R$.
Following \cite{ARHE2} (or \cite{TR} under the name filter regular element)
we call an element 
$x \in R_{(1,0)}$ \it an almost regular element for $R$ \rm
(with respect to the $x$-degree) if 
$$
(0:_R x)_{(a,\ast)}=0 \text{ for } a \gg 0.
$$

A sequence $x_1,\ldots,x_t \in R_{(1,0)}$ is \it an almost regular sequence \rm
(with respect to the $x$-degree) if for all $i \in [t]$
the $x_i$ is almost regular for $R/(x_1,\ldots,x_{i-1})R$.

Analogue
we call an element 
$y \in R_{(0,1)}$ \it an almost regular element for $R$ \rm
(with respect to the $y$-degree) if 
$$
(0:_R y)_{(\ast,b)}=0 \text{ for } b \gg 0.
$$

A sequence $y_1,\ldots,y_t \in R_{(0,1)}$ is \it an almost regular sequence \rm 
(with respect to the $y$-degree) if for all $i \in [t]$
the $y_i$ is almost regular for $R/(y_1,\ldots,y_{i-1})R$.

It is well-known that, provided $|K|=\infty$, after a generic choice of coordinates
we can achieve that 
a $K$-basis of $R_{(1,0)}$ is almost regular for $R$ 
with respect to the $x$-degree
and
a $K$-basis of $R_{(0,1)}$ is almost regular for $R$ 
with respect to the $y$-degree.
For the convenience of the reader we give a proof of this fact, which
follows from the following lemma (see also \cite{TR}).

\begin{Lemma}
Let $R$ be a bigraded algebra. If $\dim_K R_{(1,0)} > 0$ $(\dim_K R_{(0,1)} > 0)$, 
then there exists an element $x \in R_{(1,0)}$ $(y \in R_{(0,1)})$ which
is almost regular for $R$.
\end{Lemma}
\begin{proof}
By symmetry it is enough to prove the existence of $x$. 
We claim that it is possible to choose $0 \neq x \in R_{(1,0)}$ such that
for all $Q \in \Ass_S(0:_R x)$ one has $Q \supseteq \mm_x$.
It follows that $\Rad_S(\Ann_S(0:_R x))\supseteq \mm_x$. Hence 
there exists an integer $i$ such that $\mm_x^i(0:_R x)=0$ and this proves the lemma.

It remains to show the claim. 
If $P \supseteq \mm_x$ for all $P \in \Ass_S(R)$, 
then we may choose $0 \neq x \in R_{(1,0)}$ arbitrary 
because $\Ass_S(0:_R x)\subseteq\Ass_S(R)$. 
Otherwise there exists an ideal $P \in \Ass_S(R)$ with $P \not\supseteq \mm_x$.
In this case we may choose $x \in R_{(1,0)}$ such that
$$
x \not\in \bigcup_{P \in \Ass_S(R), P \not\supseteq \mm_x} P
$$
since $|K|=\infty$.
Let $Q \in \Ass_S(0:_R x)$ be arbitrary. Then
$x \in Q$ because $x \in \Ann_S(0:_R x)$.
We also have that $Q \in \Ass_S(R)$
and this implies that $Q \supseteq \mm_x$
by the choice of $x$.
This gives the claim.
\end{proof}

Let $\xb$ and $\yb$ be almost regular for $R$ 
with respect to the $x$- and $y$-degree. Define
$$
s^x_i=\max(\{a \colon (0:_{R/(x_1,\ldots,x_{i-1})R}x_i)_{(a,\ast)}\neq 0\} \cup \{0\}), 
\quad
s^x=\max\{s^x_1,\ldots,s^x_n\}
$$
and
$$
s^y_i=\max(\{b \colon (0:_{R/(y_1,\ldots,y_{i-1})R}y_i)_{(\ast,b)}\neq 0\}\cup \{0\}), 
\quad
s^y=\max\{s^y_1,\ldots,s^y_m\}.
$$

The following theorem characterizes the $x$- and $y$-regularity.
It is the analogue of the corresponding graded version in \cite{ARHE2}.

For its proof we consider
$\Tilde{H}_0(k-1,0)=(0:_{R/(x_1,\ldots,x_{k-1})R} x_k)$ for $k \in [n]$
and 
$\Tilde{H}_0(n,k-1)=(0:_{R/(\mm_x+y_1,\ldots,y_{k-1})R} y_k)$ for $k \in [m]$.
Then the beginning of the long exact Koszul sequence 
of the Koszul homology groups of $R$
for $k \in [n]$ is modified to
$$
\ldots \to H_1(k-1,0)(-1,0) \overset{x_k}{\to} H_1(k-1,0) \to H_1(k,0)
\to \Tilde{H}_0(k-1,0)(-1,0) \to 0,
$$
and for $k \in [m]$ to
$$
\ldots \to H_1(n,k-1)(0,-1) \overset{y_k}{\to} H_1(n,k-1) \to H_1(n,k)
\to \Tilde{H}_0(n,k-1)(0,-1) \to 0.
$$

Note that for $k \in [n]$ and $i \geq 1$ one has $H_i(k,0)_{(a,\ast)}=0$ for $a \gg 0$.
Similarly for $k \in [m]$ and $i \geq 1$ one has $H_i(n,k)_{(\ast,b)}=0$ for $b \gg 0$.

\begin{Theorem}
\label{yreg}
Let $R$ be a bigraded algebra, $\xb$ almost regular for $R$ 
with respect to the $x$-degree
and 
$\yb$ almost regular for $R$ 
with respect to the $y$-degree. 
Then
$$
\reg_x(R)=s^x \text{ and } \reg_y(R)=s^y.
$$  
\end{Theorem}
\begin{proof}
By symmetry it is enough to show this theorem only for $\xb$.
Let 
$$r_{(k,0)}=\max(\{a \colon  H_i(k,0)_{(a+i,\ast)}\neq 0 \text{ for } i \in [k]\} \cup \{0\})
$$ 
for $k \in [n]$ and
$$
r_{(n,k)}=\max(\{a \colon  H_i(n,k)_{(a+i,\ast)}\neq 0 \text{ for }  i \in [n+k]\} \cup \{0\})
$$ 
for $k \in [m]$. 
Then $r_{(n,m)}=\reg_x(R)$ because $H_0(n,m)=K$.
We claim that:
\begin{enumerate}
\renewcommand{\labelenumi}{(\roman{enumi})}
\item
For $k \in [n]$ one has
$
r_{(k,0)}=\max\{s_1^x,\ldots,s_k^x\}
$.
\item
For $k \in [m]$ one has
$
r_{(n,k)}=\max\{s_1^x,\ldots,s_n^x\}.
$
\end{enumerate}
This yields the theorem.
We show (i) by induction on $k \in [n]$. For $k=1$ we have the following
exact sequence
$$
0 \to H_1(1,0) \to \Tilde{H}_0(0,0)(-1,0)\to 0 
$$
which proves this case. Let $k>1$.
Since
$$
\ldots \to H_1(k,0) \to \Tilde{H}_0(k-1,0)(-1,0)\to 0, 
$$
we get $r_{(k,0)} \geq s_k^x$. 
If $r_{(k-1,0)}=0$, then $r_{(k,0)}\geq r_{(k-1,0)}$.
Assume that $r_{(k-1,0)}>0$.
There exists an integer $i$ such that 
$H_i(k-1)_{(r_{(k-1,0)}+i,\ast)}\neq 0$.
Then by
$$
\ldots \to H_{i+1}(k,0)_{(r_{(k-1,0)}+i+1,\ast)} 
\to H_{i}(k-1,0)_{(r_{(k-1,0)}+i,\ast)}
$$
$$
\to H_{i}(k-1,0)_{(r_{(k-1,0)}+i+1,\ast)} \to \ldots
$$
we have $H_{i+1}(k,0)_{(r_{(k-1,0)}+i+1,\ast)}\neq 0$ because 
$H_{i}(k-1,0)_{(r_{(k-1,0)}+i+1,\ast)}=0$.
This gives also $r_{(k,0)} \geq r_{(k-1,0)}$.
On the other hand let $a>\max\{r_{(k-1,0)},s_k^x \}$. 
If $i \geq 2$, then by
$$
\ldots \to H_i(k-1,0)_{(a+i,\ast)} \to H_{i}(k,0)_{(a+i,\ast)}
\to H_{i-1}(k-1,0)_{(a+i-1,\ast)} \to \ldots
$$
we get $H_{i}(k,0)_{(a+i,\ast)}=0$ because $H_i(k-1,0)_{(a+i,\ast)}=H_{i-1}(k-1,0)_{(a+i-1,\ast)}=0$.
Similarly $H_{1}(k,0)_{(a+1,\ast)}=0$.
Therefore we obtain that 
$r_{(k,0)}=\max\{r_{(k-1,0)},s_k^x\}=\max\{s_1^x,\ldots,s_k^x\}$ by the induction hypothesis.

We prove (ii) also by induction on $k \in \{0,\ldots,m \}$. The case $k=0$ was shown in (i),
so let $k>0$. Assume that $a>s^x$. For $i \geq 2$ one has
$$
\ldots \to H_i(n,k-1)_{(a+i,\ast)} \to H_{i}(n,k)_{(a+i,\ast)}
\to H_{i-1}(n,k-1)_{(a+i,\ast)} \to \ldots.
$$ 
Then
we get 
$H_{i}(n,k)_{(a+i,\ast)}=0$ because $H_i(n,k-1)_{(a+i,\ast)}=H_{i-1}(n,k-1)_{(a+i,\ast)}=0$.
Similarly $H_{1}(n,k)_{(a+1,\ast)}=0$ and therefore $r_{(n,k)}\leq s^x$.
If $s^x=0$, then $r_{(n,k)}= s^x$. Assume that $0 < s^x=r_{(n,k-1)}$.
There exists an integer $i$ such that $H_i(n,k-1)_{(s^x+i,\ast)}\neq 0$. Consider
$$
\ldots \to H_i(n,k-1)_{(s^x+i,\ast)} \overset{y_k}{\to} H_{i}(n,k-1)_{(s^x+i,\ast)}
\to H_{i}(n,k)_{(s^x+i,\ast)} \to \ldots.
$$
If $H_{i}(n,k)_{(s^x+i,\ast)}=0$, then $H_{i}(n,k-1)_{(s^x+i,\ast)}=y_k H_{i}(n,k-1)_{(s^x+i,\ast)}$.
This is a contradiction by Nakayamas lemma because $H_{i}(n,k-1)_{(s^x+i,\ast)}$ is a finitely 
generated $S_y$-module.
Hence $H_{i}(n,k)_{(s^x+i,\ast)} \neq 0$ and thus $r_{(n,k)}= s^x$.
\end{proof}

\section{$d$-sequences and $s$-sequences}
The concept of a $d$-sequence was introduced by Huneke \cite{HU}.
Recall that a sequence of elements $f_1,\ldots,f_r$ in a ring is 
called a \it $d$-sequence\rm, if
\begin{enumerate}
\renewcommand{\labelenumi}{(\roman{enumi})}
\item
$f_1,\ldots,f_r$ is a minimal system of generators of the ideal $I=(f_1,\ldots,f_r)$.
\item
$(f_1,\ldots,f_{i-1}):f_i \cap I =(f_1,\ldots,f_{i-1})$.
\end{enumerate}

A result in \cite{TR} motivated the following theorem. 
For a bigraded algebra $R$ let $n_x$
denote the ideal generated by the $(1,0)$-forms
of $R$
and let 
$n_y$ denote the ideal generated by the $(0,1)$-forms of $R$.

\begin{Proposition}
\label{dseq}
Let $R$ be a bigraded algebra. Then:
\begin{enumerate}
\renewcommand{\labelenumi}{(\roman{enumi})}
\item
$\reg_x(R)=0$ if and only if 
a generic minimal system 
of generators of $(1,0)$-forms for $n_x$ is a $d$-sequence.
\item
$\reg_y(R)=0$ if and only if 
a generic minimal system 
of generators of $(0,1)$-forms for $n_y$ is a $d$-sequence.
\end{enumerate}
\end{Proposition}
\begin{proof}
By symmetry we only have to prove (i). 
Without loss of generality $\xb=x_1,\ldots, x_n$ is an almost
regular sequence for $R$ with respect to the $x$-degree
because 
a generic minimal system 
of generators of $(1,0)$-forms for $n_x$ has this property.

By \ref{yreg} one has $\reg_x(R)=0$ if and only if
$s^x=0$.
By definition of $s^x$ this is equivalent to the fact that,
for all $i \in [n]$ and all $a>0$, we have
$$ 
\biggl( \frac{(x_1,\ldots,x_{i-1}):_Rx_i}{(x_1,\ldots,x_{i-1})} \biggr)_{(a,\ast)}=0.
$$
Equivalently, for all $i \in [n]$ we obtain
$
(x_1,\ldots,x_{i-1}):_R x_i \cap n_x=(x_1,\ldots,x_{i-1}).
$
This concludes the proof.
\end{proof}
If $n_{x}$ (resp. $n_{y}$) can be generated by a $d$-sequence (not necessarily generic),
then the proof of \ref{dseq} shows that $\reg_x(R)=0$ (resp. $\reg_y(R)=0$).

For an application we recall some more definitions.
Let $I=(f_1,\ldots,f_m)\subset S_x$ be a graded ideal generated in degree $d$.
Let $R(I)$ denote  
the 
Rees algebra of $I$
and let $S(I)$ denote the symmetric algebra of $I$.
It is well known that both algebras are bigraded and have a
presentation $S/J$ for a bigraded ideal $J \subset S$. For example 
$R(I)=S_x[It]\subset S_x[t]$. Define 
$$
\phi :S \to R(I), \text{ }x_i \mapsto x_i,\text{ } y_j \mapsto f_jt,
$$  
and let $J=\Ker(\phi)$. With the assumption that $I$ is generated in one degree
we have that $J$ is a bigraded ideal. Then
we will always assume that $R(I)=S/J$.
Note that then $I^j \iso(S/J)_{(\ast,j)}(-jd)$ for all $j \in \mathbb{N}$.
Similarly we may assume that $S(I)=S/J$ for a bigraded ideal $J \subset S$. 
We also consider the finitely generated $S_x$-module $S^j(I)=(S/J)_{(\ast,j)}(-jd)$,
which we call the $j^{\rm th}$ symmetric power of $I$.

For the notion of an $s$-sequence see \cite{HERETA}.
The following results were shown in \cite{HERETA} and \cite{TR}.
\begin{Corollary}
\label{nice}
Let $I=(f_1,\ldots,f_m)\subset S_x$ be a graded ideal generated in degree $d$. Then:
\begin{enumerate}
\renewcommand{\labelenumi}{(\roman{enumi})}
\item
$I$ can be generated by an $s$-sequence $($with respect to
the reverse lexicographic order$)$
if and only if $\reg_y(S(I))=0$.
\item
$I$ can be generated by a $d$-sequence if and only if $\reg_y(R(I))=0$.
\end{enumerate}
\end{Corollary}
\begin{proof}
In \cite{HERETA} and \cite{TR} it was shown that
\begin{enumerate}
\renewcommand{\labelenumi}{(\roman{enumi})}
\item
$I$ can be generated by an $s$-sequence (with respect to
the reverse lexicographic order)
if and only if 
$n_y \subseteq S(I)$ can be generated by a $d$-sequence.
\item
$I$ can be generated by a $d$-sequence if and only if 
$n_y \subseteq R(I)$ can be generated by a $d$-sequence.
\end{enumerate}
Together with \ref{dseq} these facts conclude the proof.
\end{proof}

\section{Bigeneric initial ideals}
We recall the following definitions from \cite{ARCRNE}. 
For a monomial $x^uy^v \in S$
we set 
$$
m_x(x^uy^v)=m(u)=\max\{0, i \text{ with } u_i>0 \},
$$
$$
m_y(x^uy^v)=m(v)=\max\{0, i \text{ with } v_i>0 \}.
$$
Similarly we define for $L \subseteq [n]$,
$$
m(L)= \max\{0,  i \text{ with } i \in L\}.
$$
Let $J \subset S$ be a monomial ideal. Let $G(J)$ denote the unique minimal 
system of generators of $J$.
If $G(J)=\{z_1,\ldots,z_t\}$ with $\deg(z_i)=(u^i,v^i)\in \mathbb{N}^n \times\mathbb{N}^m$, then we set  
$m_x(J)=\max\{|u^i|\}$ and $m_y(J)=\max\{|v^i|\}$. 

$J$ is called \it bistable \rm if 
for all monomials $z \in J$, all $i \leq m_x(z)$, all $j \leq m_y(z)$ one has
$x_iz/x_{m_x(z)}\in J$ and $y_jz/y_{m_y(z)}\in J$.
$J$ is called \it strongly bistable \rm if 
for all monomials $z \in J$, all $i \leq s$ with $x_s$ divides $z$, all $j \leq t$ 
with $y_t$ divides $z$ one has
$x_iz/x_{s}\in J$ and $y_jz/y_{t}\in J$.
\begin{Lemma}
\label{almost}
Let $J \subset S$ be a bistable ideal and $R=S/J$. Then:
\begin{enumerate}
\renewcommand{\labelenumi}{(\roman{enumi})}
\item
$x_n,\ldots,x_1$ is an almost regular sequence for $R$ with respect to the $x$-degree.
\item
$y_m,\ldots,y_1$ is an almost regular sequence for $R$ with respect to the $y$-degree.
\end{enumerate}
\end{Lemma}
\begin{proof}
This follows easily from the fact that $J$ is bistable.
\end{proof}

We fix a term order $>$ on $S$ by defining 
$x^uy^v>x^{u'}y^{v'}$ if either $(|u|+|v|,|v|,|u|)>(|u'|+|v'|,|v'|,|u'|)$ lexicographically or 
$(|u|+|v|,|v|,|u|)=(|u'|+|v'|,|v'|,|u'|)$ and $x^uy^v>x^{u'}y^{v'}$ reverse lexicographically
induced by $y_1>\ldots>y_m>x_1>\ldots>x_n$ 
(see \cite{EI} for details on monomial orders).
For a bigraded ideal $J$ let $\ini(J)$ denote the monomial ideal generated
by $\ini(f)$ for all $f \in J$. 
In \cite{ARCRNE} the bigeneric initial ideal $\bigin(J)$ was constructed in 
the following way: 
For $t \in \mathbb{N}$ let $\GL(t,K)$ be the general linear group
of the $t \times t$-matrices with entries in $K$.
Let $G=\GL(n,K)\times \GL(m,K)$ and $g=(d_{ij},e_{kl})\in G$.
Then $g$ defines an $S$ automorphism by extending $g(x_j)=\sum_i d_{ij}x_i$ and
$g(y_l)=\sum_k e_{kl}y_k$. There exists a non-empty Zariski open set $U \subset G$ such
that for all $g \in U$ we have $\bigin(J)=\ini(gJ)$. We call these $g \in U$ \it generic \rm for $J$. 
If char$(K)=0$, then $\bigin(J)$ is strongly bistable for every bigraded ideal $J$. 
See for example \cite{ARHE2} for similar results in the graded case.
\begin{Proposition}
\label{bigin3}
Let $J \subset S$ be a bigraded ideal. If $\char(K)=0$, then
$$
\reg_x(S/J)=\reg_x(S/\bigin(J)).
$$
\end{Proposition}
\begin{proof}
Set $\xb=x_n, \ldots, x_1$,
choose $g \in G$ generic for $J$ and let $\Tilde{\xb}=\Tilde{x}_n,\ldots,\Tilde{x}_1$ 
such that $x_i=g(\Tilde{x}_i)$.
We may assume that the sequence 
$\Tilde{\xb}$ is almost regular for $S/J$ with respect to the $x$-degree. 
Furthermore by \ref{almost} the sequence 
$\xb$ is almost regular for $S/\bigin(J)$ with respect to the $x$-degree. 
We have
$$
(0:_{S/((\Tilde{x}_n,\ldots,\Tilde{x}_{i+1})+J)}\Tilde{x_i})
\iso
(0:_{S/((x_n,\ldots,x_{i+1})+g(J))}x_i).
$$
It follows from \cite[15.12]{EI} that
$$
(0:_{S/((x_n,\ldots,x_{i+1})+g(J))}x_i)
\iso
(0:_{S/((x_n,\ldots,x_{i+1})+\bigin(J))}x_i).
$$
By \ref{yreg} we get the desired result.
\end{proof}

\begin{Remark}
\label{xbi}\rm
\begin{enumerate}
\renewcommand{\labelenumi}{(\roman{enumi})}
\item
In general it is not true that
$$
\reg_y(S/J)=\reg_y(S/\bigin(J)).
$$
For example let $S=K[x_1,\ldots,x_3,y_1,\ldots,y_3]$ and
$J=(y_2x_2-y_1x_3,y_3x_1-y_1x_3)$.
Then the minimal bigraded free resolution of $S/J$ is given by 
$$
0 \to S(-2,-2) \to S(-1,-1)\dirsum S(-1,-1) \to S \to 0.
$$
Therefore 
$
\reg_x(S/J)=0 \text{ and } \reg_y(S/J)=0.
$
On the other hand $\bigin(J)=(y_2x_1,y_1x_1,y_1^2x_2)$
with the minimal bigraded free resolution of $S/\bigin(J)$
$$
0 \to S(-2,-2)\dirsum S(-1,-2) 
$$
$$
\to S(-1,-1)\dirsum S(-1,-1)\dirsum S(-1,-2) \to S \to 0.
$$
Hence
$
\reg_x(S/\bigin(J))=0 \text{ and } \reg_y(S/\bigin(J))=1.
$
\item
It is easy to calculate the $x$- and the $y$-regularity
of bistable ideals. 
In fact, in \cite{ARCRNE} it was shown that for a bistable ideal 
$J \subset S$ we have
$$\reg_x(J)=m_x(J) \text{ and }\reg_y(J)=m_y(J).$$
\end{enumerate}
\end{Remark}

\section{Regularity of powers and symmetric powers of ideals}

Consider a bigraded algebra $R=S/J$ where $J$ is a bistable ideal.
Note that by \ref{almost} the sequence $x_n,\ldots,x_1$ is almost regular for $R$
with respect to the $x$-degree. 
For $i \in [n]$ and $j \geq 0$ we define
$$
m_j^i=\max(\{a \in\mathbb{N} \colon  (0:_{R/(x_n,\ldots,x_{i+1})R}x_i)_{(a,j)}\neq 0 \} \cup \{0 \}).
$$
Furthermore for a bistable ideal $J$ and $v \in \mathbb{N}^n$ we set 
$J_{(\ast,v)}=I_vy^v$ where $I_v \subset S_x$ is again a monomial ideal, 
which is stable
in the usual sense, that is if $x^u \in I_v$, then $x_ix^u/x_{m(u)} \in I_v$ for 
$i \leq m(u)$.

\begin{Proposition}
\label{constant}
Let $J \subset S$ be a bistable ideal and $R=S/J$. Then:
\begin{enumerate}
\renewcommand{\labelenumi}{(\roman{enumi})}
\item
For every $i \in [n]$ and for $j \geq 0$ we have $m^i_j \leq \max\{m_x(J)-1,0 \}$.
\item
For every $i \in [n]$ and for $j \geq m_y(J)$ we have $m^i_j=m^i_{m_y(J)}$.
\end{enumerate}
\end{Proposition}
\begin{proof}
If $G(J)=\{x^{u^k}y^{v^k} \colon  k=1,\ldots,r \}$, then $I_v=(x^{u^k} \colon v^k \preceq v)$
for $v \in \mathbb{N}^n$. 
This means that for all $x^u \in G(I_v)$ one has
$|u|\leq m_x(J)$. For fixed $v$ with $|v|=j$ we have 
$$
(0:_{R/(x_n,\ldots,x_{i+1})R}x_i)_{(\ast,v)}
=\frac{((x_n,\ldots,x_{i+1})+I_v:_{S_x}x_i)}{(x_n,\ldots,x_{i+1})+I_v}y^v.
$$
As a $K$-vector space
$$
\frac{((x_n,\ldots,x_{i+1})+I_v:_{S_x}x_i)}{(x_n,\ldots,x_{i+1})+I_v}y^v
=\Dirsum_{x^u \in G(I_v),m(u)=i} K(x^u/x_{m(u)}) y^v
$$
because $I_v$ is stable. Thus
$
m^i_j \leq \max\{m_x(J)-1,0 \},
$
which is (i).

To prove (ii) we replace $J$ by $J_{(\ast,\geq m_y(J))}$ and may assume
that $J$ is generated in $y$-degree $t=m_y(J)$.
Then $G(J)=\{x^{u^k}y^{v^k} \colon k=1,\ldots,r \}$ where $|v^k|=t$ for all $k \in [r]$. 
Let $|u^k|$ be maximal with $m(u^k)=i$ and define $c^i=\max\{|u^k|-1,0 \}$.
We show that $m^i_j=c^i$ for $j \geq t$ and this gives (ii).
By a similar argument as in (i) we have $m^i_{s+t}\leq c^i$ for all $s \geq 0$. If
$c^i=0$, then $m^i_{s+t}=0$. Assume that $c^i \neq 0$.
We claim that 
$$
(*) \quad 0 \neq [(x^{u^k}/x_i)y^{v^k}y_n^s] \in (0:_{R/(x_n,\ldots,x_{i+1})R}x_i)_{(\ast,s+t)}
\text{ for } s \geq 0.
$$
Assume this is not the case, then either 
$$
(x^{u^k}/x_i)y^{v^k}y_n^s=x_lx^{u'}y^{v'}
$$ 
for some $u',v'$ and $l \geq i+1$ which contradicts
to $m(u^k)=i$.
Or
$$
(x^{u^k}/x_i)y^{v^k}y_n^s=x^{u^{k'}}y^{v^{k'}}x^{u'}y^{v'}
$$
for $x^{u^{k'}}y^{v^{k'}} \in G(J)$. It follows that $|v'|=s$.
Let $k_1$ be the largest integer $l$ such that $y_n^l|y^{v^{k'}}$.
Then
$$
(x^{u^k}/x_i)y^{v^k}=((x^{u^{k'}}y^{v^{k'}}x^{u'})/y_n^{k_1})y^{v'}/y_n^{s-k_1}\in J
$$
because $J$ is bistable, and this is again a contradiction. Therefore
$(*)$ is true and we get 
$m^i_{s+t}\geq c^i$ for $s \geq 0$.
This concludes the proof.
\end{proof}
\begin{Remark}\rm
This proposition could also be formulated 
by changing the roles of $\xb$ and $\yb$. 
\end{Remark}

Let $A$ be a standard graded $K$-algebra. For a finitely generated graded 
$A$-module $M$ the usual Castelnuovo-Mumford regularity is defined to be
$$
\reg^A(M)=\sup\{r \in \mathbb{Z} \colon  \beta^A_{i,i+r}(M)\neq 0 \text{ for some integer } i \}.
$$
In \cite{CUHETR} and \cite{KO} it was shown that for a graded ideal
$I \subset S_x$ the function $\reg^{S_x}(I^j)$ is a linear function 
$pj+c$ for $j \gg 0$.
In the case that $I$ is generated in one degree
we give an upper bound for $c$ 
and find an integer $j_0$ for which $\reg^{S_x}(I^j)$
is a linear function for all $j \geq j_0$.

\begin{Theorem}
\label{power}
Let $I=(f_1,\ldots,f_m)\subset S_x$ be a graded ideal generated in degree $d \in \mathbb{N}$.
Let $R(I)=S/J$ for a bigraded ideal $J$. Then:
\begin{enumerate}
\renewcommand{\labelenumi}{(\roman{enumi})}
\item
$\reg^{S_x}(I^j) \leq jd+ \reg_x^S(R(I))$.
\item
$\reg^{S_x}(I^j)=jd+c$ for $j \geq m_y(\bigin(J))$ and some constant $0 \leq c \leq \reg_x^S(R(I))$. 
\end{enumerate}
\end{Theorem}
\begin{proof}
We choose an almost regular sequence 
$\Tilde{\xb}=\Tilde{x}_n,\ldots,\Tilde{x}_1$ for $R(I)$ over $S$ 
with respect to the $x$-degree. 
We have that for all $j \in \mathbb{N}$ the sequence $\Tilde{\xb}$ is almost regular 
for $I^j$ over $S_x$ in the sense of \cite{ARHE2}
because $R(I)_{(\ast,j)}(-dj)\iso I^j$ as graded $S_x$-modules
and
$$
(0:_{R(I)/(\Tilde{x}_{n},\ldots,\Tilde{x}_{i+1})R(I)}\Tilde{x}_i)_{(\ast,j)}(-dj)
=(0:_{I^j/(\Tilde{x}_{n},\ldots,\Tilde{x}_{i+1})I^j}\Tilde{x}_i).
$$
Define $m^i_j$ for $\bigin(J)$ as in \ref{constant}. Since
$$
(0:_{R(I)/(\Tilde{x}_{n},\ldots,\Tilde{x}_{i+1})R(I)}\Tilde{x}_i)\iso
(0:_{S/((x_n,\ldots,x_{i+1})+\bigin(J))}x_i),
$$
it follows that
$$
jd+m^i_j=r^i_j=\max(\{l \colon (0:_{I^j/(\Tilde{x}_{n},\ldots,\Tilde{x}_{i+1})I^j}\Tilde{x}_i)_l \neq 0\}
\union \{ 0 \}).
$$
By a characterization of the regularity of graded modules in \cite{ARHE2} we have
$
\reg^{S_x}(I^j)=\max\{jd,r^1_j,\ldots,r^n_j\}.
$
Hence the assertion follows from \ref{bigin3}, \ref{xbi}(ii) and \ref{constant}.
\end{proof}

Similarly as in \ref{power} one has:
\begin{Theorem}
\label{power2}
Let $I=(f_1,\ldots,f_m)\subset S_x$ be a graded ideal generated in degree $d \in \mathbb{N}$.
Let $S(I)=S/J$ for a bigraded ideal $J$. Then:
\begin{enumerate}
\renewcommand{\labelenumi}{(\roman{enumi})}
\item
$\reg^{S_x}(S^j(I)) \leq jd+ \reg_x^S(S(I))$.
\item
$\reg^{S_x}(S^j(I))=jd+c$ for $j \geq m_y(\bigin(J))$ and some constant $0 \leq c \leq \reg_x^S(S(I))$. 
\end{enumerate}
\end{Theorem}

Blum \cite{BL} proved the following with different methods. 
\begin{Corollary}
\label{power3}
Let $I=(f_1,\ldots,f_m)\subset S_x$ be a graded ideal generated in degree $d \in \mathbb{N}$.
\begin{enumerate}
\renewcommand{\labelenumi}{(\roman{enumi})}
\item
If $\reg_x(R(I))=0$, then
$
\reg^{S_x}(I^j) = jd
$
for $j \geq 1$.
\item
If $\reg_x(S(I))=0$, then
$
\reg^{S_x}(S^j(I)) = jd
$
for $j \geq 1$.
\end{enumerate}
\end{Corollary}
\begin{proof}
This follows from \ref{power} and \ref{power2}.
\end{proof}

Next we give a more theoretic bound for 
the regularity function becoming linear.
Consider a bigraded algebra $R$. Let $y$ be almost regular for all $\Tor^S_i(S/\mm_x,R)$ 
with respect to the $y$-degree.
Define 
$$w(R)=\max\{b \colon (0:_{\Tor^S_i(S/\mm_x,R)}y)_{(\ast,b)}\neq 0 \text{ for some } i \in[n]\}.$$

\begin{Lemma}
\label{1st}
Let $I=(f_1,\ldots,f_m)\subset S_x$ be a graded ideal generated in degree $d \in \mathbb{N}$.
\begin{enumerate}
\renewcommand{\labelenumi}{(\roman{enumi})}
\item
For $j > w(R(I))$ we have
$
\reg^{S_x}(I^{j+1}) \geq \reg^{S_x}(I^{j})+d.
$
\item
For $j > w(S(I))$ we have
$
\reg^{S_x}(S^{j+1}(I))\geq \reg^{S_x}(S^{j}(I))+d.
$
\end{enumerate}
\end{Lemma}
\begin{proof}
We prove the case $R = R(I)$. 
For $j >w(R)$ one has the exact sequence 
$$
0 \to \Tor^S_i(S/\mm_x,R)_{(\ast,j)} \overset{y}{\to} \Tor^S_i(S/\mm_x,R)_{(\ast,j+1)}.
$$
In \cite[3.3]{CUHETR} it was shown that
$$
\Tor^S_i(S/\mm_x,R)_{(a,j)} \iso \Tor^{S_x}_i(K,I^j)_{a+jd}
$$
and this concludes the proof.
\end{proof}

\begin{Lemma}
\label{2nd}
Let $R$ be a bigraded algebra. Then
$$
H_\lpnt(0,m)_{(\ast,j)}=0 \text{ for }j > \reg_y(R)+m.
$$
\end{Lemma}
\begin{proof}
We know that
$$
H_\lpnt(0,m)\iso \Tor^S_{\lpnt}(S/\mm_y,R)\iso H_{\lpnt}(S/\mm_y \tensor_S F_{\lpnt})
$$
where $F_{\lpnt}$ is the minimal bigraded free resolution of $R$ over $S$. 
Let 
$$F_i=\Dirsum S(-a,-b)^{\beta^S_{i,(a,b)}(R)}.$$
Then by the definition of the $y$-regularity we have $b \leq \reg_y(R)+m$ for all
$\beta^S_{i,(a,b)}(R) \neq 0$.
Thus $(S/(\yb)\tensor_S F_i)_{(\ast,j)}=0$ for $j>\reg_y(R)+m$.
The assertion follows.
\end{proof}

We get the following exact sequences.  
\begin{Corollary}
\label{3rd}
Let $I=(f_1,\ldots,f_m)\subset S_x$ be a graded ideal generated in degree $d \in \mathbb{N}$.
\begin{enumerate}
\renewcommand{\labelenumi}{(\roman{enumi})}
\item
For $j > \reg_y(R(I))+m$ we have the exact sequence
$$
0 \to I^{j-m}(-md) \to \Dirsum_m I^{j-m+1}(-(m-1)d) \to \ldots \to \Dirsum_m I^{j-1}(-d) \to I^{j} \to 0.   
$$
\item
For $j > \reg_y(S(I))+m$ we have the exact sequence
$$
0 \to S^{j-m}(I)(-md) \to \Dirsum_m S^{j-m+1}(I)(-(m-1)d) \to 
$$
$$
\ldots \to \Dirsum_m S^{j-1}(I)(-d) \to S^{j}(I) \to 0.   
$$
\end{enumerate}
\end{Corollary}
\begin{proof}
This statement follows from \ref{2nd}
and the fact that $R(I)_{(\ast,j)}(-jd)\iso I^j$ 
or $S(I)_{(\ast,j)}(-jd)\iso S^j(I)$ respectively.
\end{proof}

\begin{Corollary}
\label{fourth}
Let $I=(f_1,\ldots,f_m)\subset S_x$ be a graded ideal generated in degree $d \in \mathbb{N}$. Then:
\begin{enumerate}
\renewcommand{\labelenumi}{(\roman{enumi})}
\item
For 
$j \geq\max\{\reg_y(R(I))+m,w(R(I))+m\}$ we have
$$
\reg^{S_x}(I^{j+1})=d+\reg^{S_x}(I^{j}).
$$
\item
For 
$j \geq\max\{\reg_y(S(I))+m,w(S(I))+m\}$ we have
$$
\reg^{S_x}(S^{j+1}(I))=d+\reg^{S_x}(S^{j}(I)).
$$
\end{enumerate}
\end{Corollary}
\begin{proof}
We prove the corollary for $R(I)$. 
By \ref{3rd} and by standard arguments (see \ref{reso} for the bigraded case)
we get that for $j \geq \reg_y(R(I))+m$
$$
\reg^{S_x} (I^{j+1})\leq \max\{\reg^{S_x}(I^{j+1-i})+id-i +1 \colon i \in [m]\}.
$$
Since $j+1-i > w(R(I))$, it follows from \ref{1st} that
$$
\reg^{S_x}(I^{j+1-i})\leq \reg^{S_x}(I^{j+1-i+1})-d\leq \ldots \leq\reg^{S_x}(I^{j+1})-id. 
$$
Hence $\reg^{S_x}(I^{j+1})=\reg^{S_x}(I^{j})+d$.
\end{proof}

We now consider a special case
where $\reg^{S_x}(I^j)$ can be computed precisely.

\begin{Proposition}
\label{ci}
Let $R=S/J$ be a bigraded algebra which is a complete intersection. 
Let $\{z_1,\ldots,z_t\}$ be a homogeneous minimal system of generators of $J$
which is a regular sequence. 
Assume that 
$\deg_x(z_t)\geq\ldots\geq\deg_x(z_1)>0$ and $\deg_y(z_k)=1$ for all $k \in [t]$. 
Then for all $j \geq t$
$$
\reg^{S_x}(R_{(\ast,j+1)})=\reg^{S_x}(R_{(\ast,j)}).
$$
If in addition $\deg_x(z_k)=1$ for all $k \in [t]$,
then for $j \geq 1$
$$
\reg^{S_x}(R_{(\ast,j)})=0.
$$
\end{Proposition}
\begin{proof}
The Koszul $K_{\lpnt}(\zb)$ complex with respect to $\{z_1,\ldots,z_t \}$ provides a minimal 
bigraded free
resolution of $R$ because these elements form a regular sequence.
Observe that $(\ast,j)$ is an exact functor on complexes of bigraded modules.
Note that $K_{\lpnt}(\zb)_{(\ast,j)}$ is a complex of free $S_x$-modules
because 
$$
K_{i}(\zb)\iso\Dirsum_{\{j_1,\ldots,j_i\}\subseteq [t]}S(-\deg(z_{j_1})-\ldots-\deg(z_{j_i})),
$$
and 
$$
S(-a,-b)_{(\ast,j)}\iso\Dirsum_{|v|=j-b} S_x(-a)y^v
\text{ as graded $S_x$-modules}.
$$
Furthermore $K_{\lpnt}(\zb)_{(\ast,j)}$ is minimal by the additional assumption $\deg_x(z_k)>0$.
We have for $j \geq t$
$$
\reg^{S_x}(R_{(\ast,j)})=\max\{ \deg_x(z_t)+\ldots+\deg_x(z_{t-i+1}) -i \colon i \in [t] \}
$$
and this is independent of $j$. If in addition $\deg_x(z_k)=1$ for all $k$, then we obtain 
$$
\reg^{S_x}(R_{(\ast,j)})=0
\text{ for } j \geq 1.
$$
\end{proof}

Recall that a graded ideal $I$ is said to be of linear type,
if $R(I)=S(I)$. For example ideals generated by $d$-sequences
are of linear type.
Let $I=(f_1,\ldots,f_m)\subset S_x$ be a graded ideal,
which is Cohen-Macaulay of codim 2.
By the Hilbert-Burch theorem 
$S_x/I$ has a minimal graded free resolution
$$
0 \to \Dirsum_{i=1}^{m-1} S_x(-b_i) \overset{B}{\to} \Dirsum_{i=1}^{m} S_x(-a_i) \to S_x \to S_x/I \to 0
$$
where $B=(b_{ij})$ is a $m \times m-1$-matrix 
with $b_{ij}\in \mm$ and we may assume that the ideal $I$ is generated by the maximal minors of $B$.
The matrix $B$ is said to be the Hilbert-Burch matrix of $I$. 
If $I$ is generated in degree $d$,
then $S(I)=S/J$ where $J$ is the bigraded ideal 
$(\sum_{i=1}^m b_{ij}y_i \colon j=1,\ldots,m-1)$.

\begin{Corollary}
\label{burch}
Let $I=(f_1,\ldots,f_m)\subset S_x$ be a graded ideal generated in degree $d \in \mathbb{N}$,
which is Cohen-Macaulay of codim 2  with $m \times m-1$ Hilbert-Burch matrix $B=(b_{ij})$ and of linear type.
Then for $j \geq m-1$
$$
\reg^{S_x}(I^{j+1})=\reg^{S_x}(I^{j})+d.
$$
If additionally $\deg_x(b_{ij})=1$ for $b_{ij}\neq 0$, 
then the equality holds for $j \geq 1$.
\end{Corollary}
\begin{proof}
Since $I$ is of linear type, we have 
$R(I)=S(I)=S/J$ with the ideal
$J=(\sum_{i=1}^m b_{ij}y_i \colon j=1,\ldots,m-1)$.
One knows that (Krull-) $\dim(R(I))=n+1$. 
Since $J$ is defined by $m-1$ equations, 
we conclude that $R(I)$ is a complete intersection.
Now apply \ref{ci}.
\end{proof}

\section{Bigraded Veronese algebras}
Let $R$ be a bigraded algebra and fix $\Tilde{\Delta}=(s,t) \in \mathbb{N}^2$ with $(s,t)\neq (0,0)$.
We call 
$$
R_{\Tilde{\Delta}}=\Dirsum_{(a,b)\in \mathbb{N}^2}R_{(as,bt)}
$$ 
the 
\it bigraded Veronese algebra of $R$ \rm with respect to $\Tilde{\Delta}$
(see for example \cite{EIRETO} for the graded case and
\cite{COHETRVA} for similar constructions in the bigraded case).
Note that $R_{\Tilde{\Delta}}$ is again a bigraded algebra.
We want to relate $\reg_x^{S_{\Tilde{\Delta}}}(R_{\Tilde{\Delta}})$ 
and $\reg_y^{S_{\Tilde{\Delta}}}(R_{\Tilde{\Delta}})$
to $\reg_x^S(R)$ and $\reg_y^S(R)$. We follow the way presented in 
\cite{COHETRVA} for the case of diagonals.

\begin{Lemma}
\label{reso}
Let $R$ be a bigraded algebra and
$$
0 \to M_r \to \ldots \to M_0 \to N \to 0
$$
be an exact complex of finitely generated bigraded $R$-modules. Then
$$
\reg_x^R(N) \leq \sup \{\reg_x^R(M_k)-k \colon 0 \leq k \leq r \}
$$ 
and
$$
\reg_y^R(N) \leq \sup \{\reg_y^R(M_k)-k \colon 0 \leq k \leq r \}.
$$ 
\end{Lemma}
\begin{proof}
We prove 
by induction on $r \in \mathbb{N}$ the inequality above for $\reg_x^R(N)$.
The case $r=0$ is trivial.
Now let $r>0$, and consider
$$
0 \to N' \to M_0 \to N \to 0
$$
where $N'$ is the kernel of $M_0 \to N$.
Then for every integer $a$ we have the exact sequence
$$
 \ldots \to \Tor_i^R(M_0,K)_{(a+i,\ast)}\to \Tor_i^R(N,K)_{(a+i,\ast)}
 \to \Tor_{i-1}^R(N',K)_{(a+1+i-1,\ast)} \to \ldots
$$ 
We get
$$
\reg_x^R(N)\leq\sup\{\reg_x^R(M_0),\ \reg_x^R(N')-1\}\leq\sup \{\reg_x^R(M_k)-k \colon 0 \leq k \leq r \}
$$
where the last inequality follows from the induction hypothesis.
Analogously we obtain the inequality for $\reg_y^R(N)$.
\end{proof}

\begin{Lemma}
\label{together}
Let $A$ and $B$ be graded $K$-algebras,
$M$ be a finitely generated graded $A$-module 
and
$N$ be a finitely generated graded $B$-module.
Then
$M \tensor_K N$ is a finitely generated bigraded $A \tensor_K B$-module with
$$
\reg_x^{A \tensor_K B}(M \tensor_K N)=\reg^{A}(M) \text{ and } 
\reg_y^{A \tensor_K B}(M \tensor_K N)=\reg^{B}(N).
$$
\end{Lemma}
\begin{proof}
Let $F_{\lpnt}$ be the minimal graded free resolution of $M$ over $A$ 
and $G_{\lpnt}$ be the minimal graded free resolution of $N$ over $B$.
Then
$H_{\lpnt}=F_{\lpnt}\tensor_K G_{\lpnt}$ is the minimal bigraded free resolution
of $M \tensor_K N$ over $A \tensor_K B$ with
$H_i=\Dirsum_{k+l=i}F_k \tensor G_l$. 
Since $A(-a)\tensor_K B(-b)= (A \tensor_K B)(-a,-b)$,
the assertion follows.
\end{proof}

\begin{Theorem}
\label{reg2}
Let $R$ be a bigraded algebra, $\Tilde{\Delta}=(s,t)\in \mathbb{N}^2$ with $(s,t)\neq (0,0)$. Then
$$
\reg_x^{S_{\Tilde{\Delta}}}(R_{\Tilde{\Delta}})\leq 
\max\{c \colon c=\lceil a/s \rceil-i,  \beta^S_{i,(a,b)}(R) \neq 0 \text{ for some }i,b \in\mathbb{N}\}
$$
and
$$
\reg_y^{S_{\Tilde{\Delta}}}(R_{\Tilde{\Delta}})\leq 
\max\{c \colon c=\lceil b/t \rceil-i,  \beta^S_{i,(a,b)}(R) \neq 0 \text{ for some }i,a \in\mathbb{N}\}.
$$
\end{Theorem}
\begin{proof}
By symmetry it suffices to show the inequality for $\reg_x^{S_{\Tilde{\Delta}}}(R_{\Tilde{\Delta}})$.
Let 
$$
0 \to F_r \to \ldots \to F_0 \to R \to 0 
$$
be the minimal bigraded free resolution of $R$ over $S$. 
Since $()_{\Tilde{\Delta}}$ is an exact functor,
we obtain the exact complex of finitely generated $S_{\Tilde{\Delta}}$-modules
$$
0 \to (F_r)_{\Tilde{\Delta}} \to \ldots \to (F_0)_{\Tilde{\Delta}} \to R_{\Tilde{\Delta}} \to 0. 
$$
By \ref{reso} we have
$$
\reg_x^{S_{\Tilde{\Delta}}}(R_{\Tilde{\Delta}})\leq 
\max\{\reg_x^{S_{\Tilde{\Delta}}}((F_i)_{\Tilde{\Delta}})-i\}.
$$
Since
$$
F_i=\Dirsum_{(a,b)\in\mathbb{N}^2}S(-a,-b)^{\beta^S_{i,(a,b)}(R)},
$$
one has
$$
\reg_x^{S_{\Tilde{\Delta}}}((F_i)_{\Tilde{\Delta}})
=\max\{
\reg_x^{S_{\Tilde{\Delta}}}
(S(-a,-b)_{\Tilde{\Delta}}) \colon  \beta^S_{i,(a,b)}(R)\neq 0\}.
$$
We have to compute $\reg_x^{S_{\Tilde{\Delta}}}
(S(-a,-b)_{\Tilde{\Delta}})$.
Let $M_0,\ldots,M_{s-1}$ be the relative Veronese modules of $S_x$ and
$N_0,\ldots,N_{t-1}$ be the relative Veronese modules of $S_y$. That 
is $M_j=\Dirsum_{k \in \mathbb{N}}(S_x)_{ks+j}$ for 
$j=0,\ldots,s-1$ and 
$N_j=\Dirsum_{k \in \mathbb{N}}(S_y)_{kt+j}$ for 
$j=0,\ldots,t-1$.
Then
$$
S(-a,-b)_{\Tilde{\Delta}} 
=\Dirsum_{(k,l)\in \mathbb{N}^2} (S_x)_{ks-a}\tensor_K (S_y)_{lt-b}
=M_i (-\lceil a/s \rceil)\tensor_K N_j (-\lceil b/t \rceil)
$$
where $i \equiv-a$ mod $s$ for $0 \leq i \leq s-1$ and
$j \equiv-b$ mod $t$ for $0 \leq j \leq t-1$.

By \cite{ARBAHE} the relative Veronese modules over a polynomial ring have a linear resolution
over the Veronese algebra. Hence \ref{together} yields
$
\reg_x^{S_{\Tilde{\Delta}}}
(S(-a,-b)_{\Tilde{\Delta}})=\lceil a/s \rceil.
$
This concludes the proof.
\end{proof}

\begin{Corollary}
\label{nice2}
Let $R$ be a bigraded algebra. 
\begin{enumerate}
\renewcommand{\labelenumi}{(\roman{enumi})}
\item
For $s \gg 0, t \in \mathbb{N}$ and $\Tilde{\Delta}=(s,t)$ one has
$
\reg_x^{S_{\Tilde{\Delta}}}(R_{\Tilde{\Delta}})=0.
$
\item
For $t \gg 0, s \in \mathbb{N}$ and $\Tilde{\Delta}=(s,t)$ one has
$
\reg_y^{S_{\Tilde{\Delta}}}(R_{\Tilde{\Delta}})=0.
$
\end{enumerate}
\end{Corollary}

%\ \\
%\noindent
%\ \\
Tim R\"omer\\
FB6 Mathematik und Informatik\\
Universit\"at Essen\\
45117 Essen, Germany\\
tim.roemer@gmx.de

\begin{thebibliography}{99}
\bibitem{ARBAHE}
A. Aramova, S. Barcanescu and J. Herzog, 
On the rate of relative Veronese submodules,
\it Rev. Roumaine Math. Pures Appl. \rm 
\bf 40 \rm (1995), 
no. 3-4, 
243-251.
\bibitem{ARCRNE}
A. Aramova, K. Crona and E. De Negri, 
Bigeneric initial ideals, diagonal subalgebras and 
bigraded Hilbert functions, 
\it J. Pure Appl. Algebra \rm 
\bf 150 \rm (2000), 
no. 3, 
215-235. 
\bibitem{ARHE2}
A. Aramova and J. Herzog,
Almost regular sequences and Betti numbers, 
\it Amer. J. Math. \rm 
\bf 122 \rm (2000), 
no. 4, 
689-719.
\bibitem{BL}
S. Blum,
personal communication $(2000)$.
\bibitem{BRHE}
W. Bruns and J. Herzog, 
Cohen-Macaulay rings, 
revised edition,
\it Cambridge Studies in Advanced Mathematics \rm \bf 39\rm, 
Cambridge Univ. Press,
Cambridge,
1998.
\bibitem{COHETRVA}
A. Conca, J. Herzog, N. V. Trung and G. Valla,
Diagonal subalgebras of bigraded algebras and
embeddings of blow-ups of projective spaces,
\it Amer. J. Math. \rm
\bf 119 \rm (1997),
no. 4,
859-901. 
\bibitem{CUHETR}
S. D. Cutkosky, J. Herzog and N. V. Trung,
Asymptotic behaviour of the Castelnuovo-Mumford regularity,
\it Compositio Math.  \rm
\bf 118 \rm (1999),
no. 3,
243-261.
\bibitem{EI}
D. Eisenbud, 
Commutative algebra with a view toward algebraic geometry,
\it Graduate Texts in Mathematics \rm \bf 150\rm, 
Springer-Verlag,
New York,
1995.
\bibitem{EIRETO} 
D. Eisenbud, A. Reeves and B. Totaro,
Initial ideals, Veronese subrings, and rates of algebras,
\it Adv. Math. \rm
\bf 109 \rm (1994),
no. 2,
168-187.
\bibitem{HERETA}
J. Herzog, G. Restuccia and Z. Tang,
$s$-sequences and symmetric algebras,
preprint (2000).
\bibitem{HU}
C. Huneke,
The theory of $d$-sequences and powers of ideals,
\it Adv. in Math. \rm
\bf 46 \rm (1982),
no. 3,
249-279.
\bibitem{KO}
V. Kodiyalam,
Asymptotic behaviour of Castelnuovo-Mumford regularity,
\it Proc. Amer. Math. Soc. \rm
\bf 128 \rm (2000),
no. 2,
407-411.
\bibitem{TR}
N. V. Trung,
The Castelnuovo regularity of the Rees algebra and the associated graded ring,
\it Trans. Amer. Math. Soc. \rm
\bf 350 \rm (1998),
no. 7,
2813-2832. 
\end{thebibliography}
\end{document}